\documentclass[12pt]{article}

 \def\beq{\begin{equation}}
\def\eeq{\end{equation}}

\def\l{\lambda}
\def\L{\Lambda}

\def\o{\omega}
\def\O{\Omega}

\def\bc{\begin{center}}
\def\ec{\end{center}}

\def\vuoto{\ \hfill\hbox{\vbox{\hrule\hbox{\vrule
height5pt\kern5pt\vrule height5pt}\hrule}}\par\medskip\rm}

\begin{document}

\title{Maslov class and minimality in Calabi-Yau manifolds}

\author{Alessandro Arsie \thanks{e-mail: arsie@sissa.it} \\
S.I.S.S.A. - I.S.A.S. \\
Via Beirut 4 - 34013 Trieste, Italy}
\date{}
\maketitle

\begin{abstract}
Generalizing the construction of the Maslov class $[\mu_{\L}]$ for a Lagrangian 
embedding in a symplectic vector space, we prove that it is possible 
to give a consistent definition of the class $[\mu_{\L}]$ for any Lagrangian 
submanifold of a Calabi-Yau manifold. Moreover, extending a result 
of Morvan in symplectic vector spaces, we prove that $[\mu_{\L}]$ 
can be represented by $i_{H}\o$, where $H$ is the mean curvature 
vector field of the Lagrangian embedding and $\o$ is the Kaehler 
form associated to the Calabi-Yau metric. Finally, we conjecture a 
 generalization of the Maslov class for Lagrangian 
submanifolds of any symplectic manifold, via the mean curvature 
representation.

\end{abstract}



\section{Introduction}

The Maslov class $[\mu_{\L}]$ of a Lagrangian embedding $j:\L 
\hookrightarrow V$ in the standard Euclidean symplectic vector space $V$ 
has been constructed by Maslov in the study of global patching problem 
for asymptotic solutions of some PDEs (see \cite{M} for further 
details on this point of view). Subsequently, this cohomological class 
has found applications in the analysis of several quantization 
procedure, starting from \cite{A1} up to recent aspects on its 
relations with asymptotic, semiclassical and geometric quantization, 
for which we refer to \cite{KM}. In spite of this, there are several problems 
in the very definition of the Maslov class for Lagrangian 
submanifolds of generic symplectic manifolds.

In \cite{Mor} it has been proved that, for a Lagrangian embedding
 $j:\L \hookrightarrow V$ in a Euclidean symplectic vector space $(V,\o)$, the 
 Maslov form $\mu_{\L}$ can be represented by $\mu_{\L}=i_{H}\o$, that is by 
 the contraction of the symplectic form with the  mean curvature 
 vector field $H$ of the embedding $j$. Unfortunately, the 
 very definition of Maslov form (and related class) as exposed in 
 \cite{A1}, \cite{A2} and \cite{M}, depends on the fact that the     
 Lagrangian submanifold $\L$ is embedded in a symplectic vector 
 space, in which we have chosen a projection $\pi:V\rightarrow \L_{0}$ 
 over a {\em fixed} Lagrangian subspace $\L_{0}$; then the Maslov 
 class $[\mu_{\L}]\in H^{1}(\L,R)$ can be defined as the 
 Poincare' dual to the singular locus $Z(\L)\hookrightarrow \L$, where 
 $Z(\L):=\{ \l\in \L | {\rm rk}(\pi_{*}(\l)) < {\rm max} \}\cap H_{n-1}(\L,Z)$. 
 In the classical literature it is proved that if one changes 
 projection $\pi$, that is if one changes the reference Lagrangian subspace 
 $\L_{0}$, then the {\em Maslov class} $\mu_{\L}$ does not change, 
 while its representative changes. This is achieved using the so called 
 universal Maslov class construction on the Lagrangian Grassmannian 
 $GrL(V)$, (the homogeneous space which parametrizes Lagrangian 
 subspaces of $(V,\o)$, see \cite{A1}, \cite{A2} and \cite{KM}).
 These formulations depend heavily on the linear structure of the 
 ambient manifold $V$; in particular it is assumed that $V$ is 
 endowed with the trivial connection. Therefore,  it seems difficult 
 even to define the Maslov class for Lagrangian submanifolds of 
 symplectic manifolds, which are not vector spaces. For instance,
 it is possible to define the Maslov class of a Lagrangian embedding via 
the so called generating functions, or 
their generalization (Morse families), for which we refer to 
\cite{M}, and particularly \cite{W}. In this way, one obtains a 
notion of Maslov class for Lagrangian submanifolds embedded in 
{\em any} cotangent bundle $T^{*}M$ over a Riemannian manifold $M$,
constructing a $Z$-valued 
${\check{\rm C}}$ech cocycle, starting from the signature of the 
Hessian of a Morse family; however this construction depends strongly 
on the choice of a Òbase manifoldÓ ($M$ in the case of the cotangent 
bundle) and does not seem to be generalizable to Lagrangian embedding 
in any symplectic manifold. (See \cite{W} for more details on this 
kind of construction).

 Recently (see \cite{F}), Fukaya has shown how to define a 
 Maslov index for closed loops on Lagrangian submanifolds of a quite 
 general class of symplectic manifolds, the so called pseudo-Einstein 
 symplectic manifolds. The construction is 
 developed using non trivial assumptions on the structure of the 
 ambient manifold  and is carried on only for a particular subclass 
 of Lagrangian submanifolds; moreover, there is no explicit reference to the 
 corresponding Maslov class. 
 
 In this paper we show that, whenever the ambient manifold is 
 Calabi-Yau, it is possible to give a consistent definition of Maslov 
 class for its Lagrangian submanifolds, generalizing the approach of 
 Arnol'd with the so called universal Maslov class. In this framework, 
 we show that it is possible to generalize the result of Morvan and 
 then we comment on various consequences of our construction, in 
 particular on the possible definition of Maslov class for Lagrangian 
 embedding in any symplectic manifold.

\section{The Maslov class for Lagrangian embedding in Calabi-Yau}

Let us briefly recall the standard construction of the Maslov class 
$\mu_{\L}$, for a Lagrangian submanifold $\L$, embedded in a 
symplectic vector space $(V,\o)$, of real dimension $2n$: first of 
all, one considers the tangent spaces to $\L$ as (affine) subspaces 
of $V$. Then, using the trivial parallel displacement one transports every 
tangent plane in a fixed point $P$ of $V$, (for example the origin). 
Now, one has to consider the Lagrangian Grassmannian $GrL(T_{P}V)$, 
which by definition parametrizes all Lagrangian subspaces of $T_{P}V$. 
Using the trivial connection, we have thus obtained a map:
\[ G: \L \longrightarrow GrL(T_{P}V). \]
It is easy to see  (\cite{A1}, \cite{A2}),that $GrL(T_{P}V)$ has the 
natural structure of the homogeneous space $ \displaystyle \frac{U(n)}{O(n)}$; then by the 
standard tool of the exact homotopy sequence for a fibration 
(see\cite{BT}), it 
is proved that $\pi_{1}(GrL(T_{P}V))\cong  Z$. In fact, having 
fixed a Lagrangian plane $\L_{0}$ in $T_{P}V$, all other Lagrangian 
planes are obtained via a unitary automorphism $A\in U(n)$. Obviously, 
we have a fibration:
\[ SU(n) \longrightarrow U(n)\stackrel{det}{\longrightarrow} S^{1}, \]
but this does not descend to $GrL(T_{P}S)$, since we have to quotient out 
the possible orthogonal automorphisms. However, since the square of 
the determinant of an orthogonal automorphism is always 1, we have a 
well defined map:
\[ det^{2}: GrL(T_{P}S) \longrightarrow S^{1}, \]
which sits in the following commutative diagram of fibrations:
\bc
\[ \begin{array}{cccccc}
  SO(n) & \longrightarrow & O(n) & 
  \stackrel{det}{\longrightarrow}& S^{0} &  \\
 \downarrow &     &   \downarrow &     &  \downarrow &    \\
 SU(n) &  \longrightarrow  &  U(n) & 
 \stackrel{det}{\longrightarrow}&  S^{1} &   \\
 \downarrow &     &  \downarrow  &     & \downarrow &  z^{2}  \\
 GrSL( C^{n}) &\longrightarrow & GrL(C^{n}) & 
 \stackrel{det^{2}}{\longrightarrow} & S^{1} &   \\
 \end{array}  \]
 \ec

In this diagram the space $GrSL(C^{n})$ denotes the 
Grassmannian of special Lagrangian planes in $C^{n}$, that is 
the Grassmannian of Lagrangian planes which are {\em calibrated} by 
the top holomorphic form of $C^{n}$; the corresponding 
Lagrangian submanifolds are called special Lagrangian (see \cite{HL} 
for more details). Notice that this space is always simply connected.

Finally, using 
Hurewicz isomorphism and taking a 
generator belonging to  $H^{1}(GrL(T_{P}V), Z)$, which is thought  
 as the pull-back via $det^{2}$ of the generator  
 $[\alpha]\in H^{1}(S^{1},Z)$, one defines the 
Maslov class $[ \mu_{\L}]:=G^{*}(det^{2})^{*}[\alpha]$. Obviously, this 
construction is indipendent on the choice of the point $P$, since if 
another point is chosen it is possible to construct a homotopy in such 
a way to prove the invariance of $[\mu_{\L}]$. It is clear that, in 
this framework, the existence of the trivial connection is an (almost!) essential 
 requirement for the construction to work. In fact, we will see in this section that, to 
have a consistent definition of Maslov class it is not necessary that 
the ambient manifold is endowed with the trivial connection, but is 
sufficient that the global holonomy of the symplectic manifold is ÒsmallÓ
in a suitable sense.

From now on we restrict our attention to Lagrangian submanifolds of 
Calabi-Yau manifolds. Recall that Calabi-Yau manifolds can be defined as 
compact  Kaehler manifolds with vanishing first Chern class; recall 
also that a celebrated theorem by Yau (proving a previous conjecture 
by Calabi) implies that for every choice of the Kaehler class on a 
Calabi-Yau, there exists a unique Ricci-flat Kaehler metric. 
Moreover, while the holonomy of a Kaehler manifold is contained in 
U(n), if $g$ is the Ricci-flat metric of an n-dimensional Calabi-Yau, 
then the corresponding holonomy group is contained in SU(n). Finally, 
let us recall that, on every Kaehler manifold $(X,g,J)$ (where $g$ 
is a Kaehler metric and $J$ the integrable almost complex structure) 
the corresponding symplect or Kaehler form $\o$ is related to $g$ 
via:
\beq
\label{simplet}
\o(X,Y):=g(X,JY) \quad \forall X,Y \in \Gamma(TX), 
\eeq
and that the almost complex structure tensor $J$ is covariantly 
constant with respect to the Levi-Civita connection induced by $g$. 
Considering a Kaehler metric $g$ on a Calabi-Yau, we will always mean 
the Ricci-flat metric. Typical examples of Calabi-Yau are given by the 
zero locus of a homogeneous polynomial of degree $n+1$ in 
$ P^{n}(C)$ (whenever this locus is smooth); however it is by no means true that all 
Calabi-Yau are algebraic. For further details on this class of 
manifolds see for example \cite{Besse} and \cite{V}.

The construction of Fukaya for defining the Maslov index of closed 
loops goes as follow (see \cite{F} for details and motivations). He 
considers symplectic manifolds $(X,\o)$ which are Òpseudo-EinsteinÓ in 
the sense that there exists an integer $N$ such that $N\o=c_{1}(X)$. 
By this relation, the line bundle  $det(TX)$ is flat when restricted to every 
Lagrangian submanifold $\L$ of $X$, but Fukaya restricts further 
the class of Lagrangian submanifolds considering only the so called 
Bohr-Sommerfeld orbit $\L$ (BS-orbit for short), which are defined as 
the Lagrangian submanifolds for which the restriction of $det(TX)$ is 
not only flat, but even trivial. This implies that if we consider a 
closed loop $h:S^{1} \longrightarrow \L$ ($\L$ is a BS-orbit), then 
the monodromy  $M$ of the tangent bundle $TX$ along $h(S^{1})$ is 
contained in $SU(n)$. Then the idea is to take a path in $SU(n)$ 
joining $M$ with the identity, in order to get an induced 
trivialization of $h^{*}(TX_{|h(S^{1})})\cong S^{1}\times { 
C}^{n}$. In this trivial bundle there is a family of Lagrangian vector 
subspaces $T_{h(t)}\L$ and in this way we get a loop in $GrL(C^{n})$ and hence a 
well-defined integer (the Maslov index) $m(h)$. 
Obviously $m(h)$ is independent of the choice of the path in $SU(n)$ 
which joins $M$ to the unit, since $\pi_{1}(SU(n))\cong 1$.

Now we come to our construction. Consider embedded Lagrangian 
submanifolds $\L$ of a Calabi-Yau $(X,\o,g,J)$, where $\o,g,J$ are 
related by (\ref{simplet}). Define the {\em Lagrangian 
Grassmannization} $GrL(X)$ of $TX$ as the fibre bundle over $X$ 
obtained substituting $T_{x}X$ with $GrL(T_{x}X)$, thus:
\[ GrL(X):=\coprod_{x\in X}^{}GrL(T_{x}X) \]
and in particular:
\[ GrL(X)_{\L}:=\coprod_{x\in \L}^{}GrL(T_{x}X). \]
 Let $G(j)$ be the Gauss map, which takes $x\in \L$ in $T_{x}\L$ 
 thought as a Lagrangian subspace of $T_{x}X$. Via $G(j)$, the 
 embedding $j:\L \hookrightarrow X$ lifts to a section $G(j):\L 
 \rightarrow GrL(X)_{\L}$. We would like to define the Maslov class 
 of $\L$ via a map ${\cal M}:\L \longrightarrow S^{1}$ in the following 
 way: to every point $x\in \L$, we consider $G(j)(x)$ and then 
 through the isomorphism $ \displaystyle GrL(T_{x}X)\cong \frac{U(n)}{O(n)}$, taking the map 
 $det^{2}$ we get a point in $S^{1}$. However, as we have seen, to establish an 
 isomorphism to every space $GrL(T_{x}X)$ ($x\in \L$) with 
 $\displaystyle \frac{U(n)}{O(n)}$ we need a reference Lagrangian plane in 
 $GrL(T_{x}X)$ $\forall x \in \L$, that is we need {\em another 
 section} of $GrL(X)_{\L}$, besides $G(j)(\L)$. 
 
 To this aim, fix a point $p\in \L$, consider $T_{p}\L$ and use the 
 parallel displacement, induced by the Levi Civita connection of $g$, 
 along a system $\gamma$ of paths on $\L$ starting from $p$, 
 to construct a  reference distribution of Lagrangian planes 
 ${\cal D}_{\gamma}$ over $\L$, that is 
 another section of $GrL(X)_{\L}$. This is indeed possible, since the 
 holonomy is contained in $U(n)$, the parallel displacement is an 
 isometry for $g$ and $J$ is covariantly constant: these facts, combined with the relation 
 \ref{simplet} imply that parallel transport sends Lagrangian planes 
 in Lagrangian planes. Obviously this distribution ${\cal D}_{\gamma}$ 
 is not uniquely determined, since it depends on the choice of the 
 system of paths $\gamma$ starting from $p$. In spite of this, due to 
 the fact that the holonomy of a Calabi-Yau  metric is very 
 constrained, this dependence does not prevent us to reach our goal. 
 Indeed, consider $q\in \L$ and compare the two Lagrangian planes 
 $({\cal D}_{\gamma})_{q}$ and $({\cal D}_{\delta})_{q}$ obtained by 
 parallel transport of $T_{p}\L$ along two different paths $\gamma$ 
 and $\delta$. By the holonomy property of a Calabi-Yau metric we 
 have:
 \[ ({\cal D}_{\gamma})_{q}=M ({\cal D}_{\delta})_{q} \quad M\in SU(n). 
 \]
Thus, if $A\in U(n)$ is such that $T_{q}\L=A({\cal D}_{\gamma})_{q}$, 
then $T_{q}\L=AM({\cal D}_{\delta})_{q}$; so to every $q\in \L$ we can 
associate $A_{q}$ such that 
$G(j)(q)=T_{q}\L=A_{q}({\cal D}_{\gamma})_{q}$, where $A_{q}$ is 
determined up to multiplication by a matrix $M\in SU(n)$. At this 
point the key observation is that $det^{2}(A_{q})\in S^{1}$ is a well 
defined point, which is not affected by the ambiguity of $A_{q}$.
In this way we have a well-defined map, the {\em Maslov map}:
\[
\begin{array}{cccc}
{\cal M}: & \L & \longrightarrow & S^{1} \\

           & q & \mapsto & det^{2}(A_{q}) \\
\end{array}
\]
Take the generator $[\alpha]$ of $H^{1}(S^{1}, Z)$ represented 
by the form $\alpha:=\frac{1}{2\pi}d\theta$. Observe that the target space of 
the Maslov map, is not only topologically a circle, but even a Lie 
group, the group $U(1)$: this implies that the choice of the form 
$\frac{1}{2\pi}d \theta$ is compulsory, since it is the unique 
normalized invariant 1-form.
Now we can give the following:

{\em Definition}: Using the previous notations, we define the {\em 
Maslov form} of the Lagrangian embedding $j:\L\hookrightarrow X$ as
 $\mu_{\L}:={\cal M}^{*}\alpha$ and the corresponding {\em Maslov 
 class} as $[\mu_{\L}]={\cal M}^{*}[\alpha] \in H^{1}(\L,Z)$.

{\em Remark} 1 : The Maslov map ${\cal M}$ has been built up fixing a 
reference point $p$, from which we constructed ${\cal D}_{\gamma}$; in 
this way the map ${\cal M}$ associates to $p$ $1\in S^{1}$. It is 
clear that if one takes a different reference point $p'$, then the 
map $\cal{M}$ changes (this time $p'$ goes to 1), but the Maslov 
class and the Maslov form do not change, as it is immediate to see. In 
particular, the invariance of the Maslov form is due to the invariance 
of $\alpha$ under the action of the Lie group $U(1)$.

{\em Remark} 2 : In \cite{Trof}, Trofimov costructed a generalized 
Maslov class, as a cohomological class defined on the space of paths 
$[X,\L]$; these paths start from a fix point $x_{0}$ in a symplectic manifold $X$ and 
end to a fixed Lagrangian submanifold $\L$ of $X$. We argue that the
 the Maslov class we have just defined can be obtained as a finite dimensional 
reduction of the class built up in \cite{Trof}, when one uses the 
Levi-Civita connection induced by the Calabi-Yau metric. In fact, 
Trofimov did not use metric connections, but instead affine torsion 
free connections, preserving the symplectic structure, which are 
generally not induced by a metric.

\section{Representation of the Maslov class via the mean curvature vector  
field}

In this section, generalizing what has been proved by Morvan in 
\cite{Mor} for Lagrangian embeddings in Euclidean symplectic vector 
space, we prove the following:

{\bf Theorem}: {\em Let} $j:\L \hookrightarrow X$ {\em be a Lagrangian embedding 
in a Calabi-Yau X and let} $H\in \Gamma (N\L)$ {\em be the mean curvature 
vector field of the embedding }$j$ {\em (with repect to the Calabi-Yau 
metric),  then}:
\[ \mu_{\L}= \frac{1}{\pi}i_{H}\o, \] 
{\em where } $\o$ {\em is the Kaehler form constructed from the 
Calabi-Yau metric} $g$, and $\mu_{\L}$ {\em is the Maslov form previously defined}.

Before proving the theorem we need various preliminary results, which 
we are going to state and prove, and we need also to decompose into 
simpler pieces the action of ${\cal M}^{*}$ on $[\alpha]$.

Recall that given an embedding $j$, the associated second fundamental 
form $\sigma: T\L\times T\L \rightarrow N\L$ is a symmetric tensor 
defined by:
\[ \sigma (X,Y):=\nabla^{g}_{X}Y-\nabla^{j^{*}g}_{X}Y, \quad \forall 
X,Y\in \Gamma(T\L), \]
where $\nabla^{g}$ is the Levi-Civita connection in the ambient 
manifold, while $\nabla^{j^{*}g}$ is the connection induced on $\L$ 
via the pulled-back metric. If $\sigma$ is identically vanishing, 
then the submanifold is called totally geodesic. Taking the trace of 
$\sigma$ we get a field of normal vectors, that is the {\em mean curvature 
vector field } $H$ of the embedding $j$. Those embeddings for which 
$H$ is identically vanishing are called minimal. 

First of all we need to understand the local structure of 
$TGrL(T_{x}X)$. Fix a point $q\in \L$ and set $V:=T_{q}X$ for short.
We can prove the following:

{\bf Lemma 1}: {\em The space} $T_{\pi}GrL(V)$ {\em over a Lagrangian 
n-plane $\pi$ of $V$ can be identified with the subspace of linear 
maps $\psi:\pi \rightarrow \pi^{\bot}$ ($\pi^{\bot}$ denotes the 
orthogonal subspace in $V$ with respect to the metric $g$ in $q$) such that}:
\[ g(\psi(X),JY)=g(\psi(Y),JX), \quad \forall X,Y \in \pi. \]

Proof: First of all, we have $T_{\pi}GrL(V)\equiv S(\pi)$, where 
$S(\pi)$ is the space of all symmetric bilinear forms on $\pi$. In 
fact every $v\in T_{\pi}GrL(V)$ can be represented as 
$\displaystyle \frac{d}{dt}B(t)\pi_{|t=0}$, where $B(t)$ is a path of linear 
symplectic transformation of $V$, with the condition $B(0)=id_{V}$. 
To $v\in T_{\pi}GrL(V)$ we can associate a form $S_{v}$ given by:
\[ S_{v}(X,Y):=\o(\frac{d}{dt}B(t)X_{|t=0},Y). \]
This form is clearly bilinear and is symmetric:
\[ S_{v}(X,Y)=\o(\frac{d}{dt}B(t)X_{|t=0},B(t)Y_{|t=0})= \]
\[= \frac{d}{dt}\o(B(t)X,B(t)Y)_{|t=0}-\o(B(t)X_{|t=0}, 
\frac{d}{dt}B(t)Y_{|t=0})=0-\o(X,\frac{d}{dt}B(t)Y_{|t=0})= \]
\[ =\o(\frac{d}{dt}B(t)Y_{|t=0},X)=S_{v}(Y,X), \]
by the fact that $B(t)$ is a symplectic linear transformation of $V$ 
and by skewsymmetry of $\o$. It is easy to verify that the 
corresponding map $T_{\pi}GrL(V) \rightarrow S(\pi)$ is an isomorphism.
Moreover  we have:
\[ 
S_{v}(X,Y)=\o(\frac{d}{dt}B(t)X_{|t=0},Y)\stackrel{(\ref{simplet})}{=}g(\frac{d}{dt}B(t)X_{|t=0}, 
JY) \]
and thus, identifying $\psi:\pi \rightarrow \pi^{\bot}$ with 
$ \displaystyle \frac{d}{dt}B(t) \pi_{|t=0}$ we get the result. \vuoto
 
 By Lemma 1 it is clear that $J$ itself, restricted to $q$, can be 
 considered not only as an element of $T_{\pi}GrL(V)$ but even as an 
 invariant vector field on $GrL(V)$, that is $J_{q}\in \Gamma 
 (TGrL(V))$. Let ${e_{1},\ldots,e_{n}}$ be an orthonormal basis of 
 $\pi$ and ${f^{1},\ldots,f^{n}}$ the corresponding dual basis, in 
 such a way that ${Je_{1},\ldots,Je_{n}}$ is a basis of $\pi^{\bot}$ 
 and ${-Jf^{1},\ldots, -Jf^{n}}$ the associated dual basis. Then $J$ 
 as a vector belonging to $T_{\pi}GrL(V)$, can be represented as a 
 section of $\pi^{*}\otimes \pi^{\bot}$, that is $J=f^{i}\otimes Je_{i}$
 (Einstein summation convention is intended). From $J$ in this 
 representation one can construct a 1-form $\tilde{J}\in 
 \O^{1}(GrL(V))$ using the paring induced by the metric, that is 
 $\tilde{J}=e_{i}\otimes -Jf^{i} $. This 1-form has a quite 
 outstanding role:
 
 {\bf Lemma 2}: {\em Fix an arbitrary Lagrangian plane in $V$ in order 
 to have a map $det^{2}:GrL(V) \rightarrow S^{1}$. Then:
 \[ (det^{2})^{*}(\alpha)=\frac{1}{\pi} \tilde{J}, \]
 so that $\tilde{J}$ defines a closed form on $GrL(V)$.}
 
 Proof: It is sufficient to prove that for every $X\in T_{\pi}GrL(V)$ 
 one has $(det^{2})^{*}(\alpha)(X)=\frac{1}{\pi} \tilde{J}(X)$.
 Indeed:
 \[ (det^{2})^{*}(\alpha)(X)=(\alpha)(det^{2}_{*}(X)), \]
 so we are led to compute the tangent map to $det^{2}$. Assume for 
 simplicity that $\pi$ is the reference Lagrangian plane in the 
 isomorphism $\displaystyle GrL(V)\cong \frac{U(n)}{O(n)}$, so that it is 
 represented by the identity matrix. Then, since 
 $ \displaystyle T_{\pi}GrL(V)\cong \frac{u(n)}{o(n)}$, consider a path $\gamma: 
 (-\epsilon, \epsilon)\rightarrow u(n)$, such that $\gamma (0)= 
  O$ and such that its image in $u(n)$ has empty intersection with 
 $o(n)$ (except for the zero matrix). The exponential mapping 
 determines in this way a path in $GrL(V)$ through $\pi$. Now, we have:
 \[ 
  \frac{d}{dt}det^{2}(e^{\gamma(t)})_{|t=0}=\frac{d}{dt}det(e^{2\gamma(t)
 })_{|t=0}= \frac{d}{dt} (e^{2Tr(\gamma(t))})_{|t=0}= \]
 \[ 2Tr(\dot{\gamma}(0))=2Tr(X), \]
 where $\dot{\gamma}(0)$ is identified with the tangent vector $X$ in 
 $T_{\pi}GrL(V)$.
 Hence one gets:
 \[(det^{2})^{*}(\alpha)(X)=(\alpha)(det^{2}_{*}(X))=(\alpha)(2Tr(X))=
 \frac{1}{\pi}Tr(X). \]
 
 On the other hand, $X\in \Gamma(\pi^{*}\otimes\pi^{\bot})$, so that 
 it can be represented as $X=X_{k}^{l}f^{k}\otimes Je_{l}$; thus one 
 gets:
 \[ \tilde{J}(X)=(e_{i}\otimes -Jf^{i})(X_{k}^{l}f^{k}\otimes Je_{l})=
 X_{i}^{i}=Tr(X).  \] \vuoto
 
 Till now we have worked only locally, having fixed a point $q\in 
 \L$. To proceed we need to globalize the properties stated in lemma 1 
 and 2. Let us define the {\em 
 vertical tangent bundle} $VT(GrL(X)_{\L})$ ($VT(GrL)$ for short) over $GrL(X)_{\L}$ as:
 \[  VT(GrL(X)_{\L}):=\coprod_{x\in \L}^{} TGrL(T_{x}X); \]
 notice that this is not the tangent bundle of $GrL(X)_{\L}$, since 
 it is obtained taking the tangent bundle of the fibre only (thus the 
 name vertical). Analogously, one can define the {\em vertical 
 cotangent bundle} over $GrL(X)_{\L}$ as:
 \[ VT^{*}(GrL(X)_{\L}):=\coprod_{x\in \L}^{}T^{*}GrL(T_{x}X), \]
 (from now on denoted as $VT^{*}(GrL)$ for short).
 
 Now, by the previous reasoning  and since $J$ is covariantly constant 
 on a Kaehler manifold $X$, we have that $J$ defines a section of 
 $VT(GrL)$ and analogously $\tilde{J}$ induces a section of 
 $VT^{*}(GrL)$. In order to globalize the result of lemma 2, observe 
 that the section ${\cal D}_{\gamma}$ of 
 $GrL(X)_{\L}$ over $\L$, defined in the previous section, enables one 
 to give a well-defined map $Det^{2}:GrL(X)_{\L}\rightarrow S^{1}$ 
 (one takes as a reference Lagrangian plane in $GrL(T_{x}X)$ the 
 subspace $({\cal D}_{\gamma})_{x}$). It is clear that one gets 
 immediately  the following:
 
 {\bf Corollary 1}: {\em Under the previous notations and considering 
 the fibration $Det^{2}: GrL(X)_{\L}\rightarrow S^{1}$ induced by the 
 reference distribution ${\cal D}_{\gamma}$ one has:}
 \[ (Det^{2})^{*}(\alpha)=\frac{1}{\pi}\tilde{J} \]
 {\em where $\tilde{J}$ is viewed as a section of $VT^{*}(GrL)$.}

 Via the Gauss map we can pull-back $VT(GrL)$ to $\L$:
 \[
 \begin{array}{ccc}
 G(j)^{*}(VT(GrL)) &  & VT(GrL) \\
 \downarrow &    &\downarrow pr_{VT} \\
 \L & \rightarrow & GrL(X)_{\L} \\
 \end{array}
  \]

  {\bf Lemma 3}: {\em The bundle $G(j)^{*}VT(GrL)$ can be identified 
  with the subspace of $T^{*}\L\otimes N\L$ consisting of those 
  sections $\psi \in \Gamma (T^{*}\L\otimes N\L)$ (that is 
  $N\L$-valued 1-forms on $\L$) such that:}
  \[ g(\psi(X), JY)=g(\psi(Y), JX), \quad \forall X,Y \in 
  \Gamma(T\L). \]
  
  Proof: By the very definition of pulled-back bundle, we have that:
  \[ G(j)^{*}VT(GrL)\cong \{ (x;x',\pi, X)\in \L\times VT(GrL): \quad
  (x,T_{x}\L)=G(j)(x)= \] 
  \[ =pr_{VT}(x',\pi,X)=(x',\pi) \}, \]
  which clearly implies the constraint $x=x'$ and $T_{x}\L=\pi$ so that:
  \[ G(j)^{*}VT(GrL)\cong \coprod_{x\in \L}^{}T_{\pi=T_{x}\L}GrL(T_{x}X). \]
  On the other hand, by lemma 1:
  \[ T_{\pi=T_{x}\L}GrL(T_{x}X)\cong\{ \psi\in 
  \Gamma(T^{*}_{x}\L\otimes N_{x}\L) \quad {\rm such \quad that}: \] 
  \[ g(\psi(X), JY)=g(\psi(Y),JX), 
  \quad \forall X,Y\in T_{x}\L \},\]
  so one gets immediately the thesis. \vuoto
  
  The tangent application to the Gauss map is related to the second 
  fundamental form as shown in the following:
  
  {\bf Lemma 4}: {\em The tangent map to $G(j)$ in a point $x\in \L$ can 
  be identified with the second fundamental form $\sigma$, thought of
  as an application with values in $T^{*}\L\otimes N\L$; more exactly
  $\sigma$ takes values in the subspace $G(j)^{*}(VT(GrL))$ of 
  $T^{*}\L\otimes N\L$, in the sense that it satisfies $g(\sigma 
  (X,Y),JZ)=g(\sigma(X,Z),JY)$}.
 
 Proof: First of all, the identity $g(\sigma(X,Y),JZ)=g(\sigma(X,Z),JY)$ 
 is a consequence of the fact that Lagrangian submanifolds of K\"ahler 
 manifolds are always anti-invariant (also called totally real) 
 submanifolds of top dimension (see \cite{YK} page 35). Hence, always 
 by result of \cite{YK}, page 43, we have the desired relation. 
 Finally, the fact that the tangent map to the Gauss map can be 
 identified with the second fundamental form, via the action of the 
 almost complex structure $J$ and the metric $g$, is a classically 
 known result which can be found, for example in \cite{BC}, page 196.  \vuoto

Observe that by lemma 3 and 4, the second fundamental form 
$\sigma(X,.)$, considered as a map taking values in $T^{*}\L\otimes 
N\L$ is an element of $G(j)^{*}(VT(GrL))$. 
Let us summarize the situation in the following diagram:

\[ 
\begin{array}{ccccc}
T\L & \stackrel{G(j)_{*}}{\rightarrow} & G(j)^{*}(VT(GrL)) & 
\subset  T^{*}\L\otimes N\L &  VT(GrL) \\
\downarrow &  &  & &  \downarrow \\
\L & \stackrel{G(j)}{\rightarrow} & G(j)(\L) &  \hookrightarrow & 
GrL(X)_{\L} \\
\end{array}
\]

Denote again with $\tilde{J}$ the restriction of $\tilde{J}$ to the bundle
$G(j)^{*}(VT^{*}(GrL))$. By the previous diagram we can pull-back 
$\tilde{J}$ to a closed 1-form on $\L$ via $G(j)^{*}$:
\beq
\label{eq1}
(G(j)^{*}(\tilde{J}))(X)=\tilde{J}(G(j)_{*}(X))=\tilde{J}(\sigma(X,.)) 
\quad \forall  X \in \Gamma(T\L),
\eeq
where the last equality in equation (\ref{eq1}) is due to lemma 4 and 
the pairing between $\tilde{J}$ and $\sigma(X,.)$ is induced by the 
natural pairing between $G(j)^{*}(VT^{*}(GrL))$ and 
$G(j)^{*}(VT(GrL))$, respectively.

\bigskip

Proof of the theorem: First of all, notice that the Maslov map ${\cal 
M}:\L \rightarrow S^{1}$ can be decomposed as ${\cal M}=Det^{2}\circ 
G(j)$, as is immediate to see. Then $\mu_{\L}:={\cal 
M}^{*}(\alpha)=G(j)^{*}\circ (Det^{2})^{*}(\alpha)$ and so 
$\mu_{\L}=\frac{1}{\pi}G(j)^{*}(\tilde{J})$, by lemma 2. Now 
$\tilde{J}=e_{l}\otimes -Jf^{l}$ and $\sigma(X,.)$ can be represented 
as $\Gamma(T^{*}\L\otimes N\L)\ni 
\sigma(X,.)=\sigma_{i}^{k}(X)f^{i}\otimes Je_{k}$. In this way we 
have that for all $X \in \Gamma (T\L)$:
\[ (G(j)^{*}(\tilde{J}))(X)=(e_{l}\otimes -Jf^{l})(\sigma_{i}^{k}(X) 
f^{i}\otimes Je_{k})=\sigma_{i}^{i}(X)= \]
\[ =\sum_{i}^{}g(\sigma(X,e_{i}),Je_{i})= {\rm  \quad (by \quad lemma \quad 
4)} \quad =\sum_{i}^{}g(\sigma(e_{i},e_{i}),JX)= \]
\[ =g(H,JX) =\o(H,X)=i_{H}\o(X). \]
Hence, one gets the result:
\beq
\label{eq2}
\mu_{\L}=G(j)^{*}(\frac{1}{\pi} \tilde{J})=\frac{1}{\pi}i_{H}\o \quad 
\in H^{1}(\L, Z).
\eeq 
\vuoto

By the result of the theorem, one can give the following:

{\em Definition}: Let $\L\hookrightarrow X$ a Lagrangian embedding in 
a Calabi-Yau $X$; then the Maslov  index  $m$ of a closed 
loop $\gamma$ on $\L$ is given by:
\[ m(\gamma):=\frac{1}{\pi}\int_{\gamma}i_{H}\o \quad \in 
 Z. \]

\section{Conclusions}

Calabi-Yau manifolds have received great attention as target spaces 
for superstring compactifications. Moreover their Lagrangian and 
special Lagrangian submanifolds are now considered as the cornerstones
 for understanding the mirror symmetry phenomenon between pairs of 
Calabi-Yau spaces, both from a categorical point of view (\cite{K}), 
and from a physical-geometrical standpoint (\cite{SYZ}). Let us recall 
that special Lagrangian submanifolds $\L$ of a Calabi-Yau $X$ are exactly what 
are called BPS states or supersymmetric cycles 
 in the physical literature; on the other hand, it is known 
that special Lagrangian submanifolds are nothing else that {\em 
minimal} Lagrangian submanifolds (compare \cite{HL} page 96, where 
this is proved for special Lagrangian submanifolds of $ C^{n}$). 
From our result it turns out that the Maslov class of special 
Lagrangian submanifolds is identically vanishing; on the other hand, 
this can be seen just by considering the Grassmannian of special 
Lagrangian planes, which turns out to be diffeomorphic to 
$\displaystyle \frac{SU(n)}{SO(n)}$, hence simply connected (notice that the 
Grassmannian of special Lagrangian planes is isomorphic to the fibre 
in the fibration $det^{2}: GrL(C^{n})\rightarrow S^{1}$). 
It is then clear that the Maslov index is identically vanishing for 
all special Lagrangian submanifolds $\L$ of a Calabi-Yau $X$. We 
believe that this simple observation can enhance our understanding  of
the structure of the $A^{\infty}$-Fukaya category, 
whenever its objects are restricted to minimal Lagrangian 
submanifolds (see \cite{F} for a definition of $A^{\infty}$ category, 
and \cite{K} for its application in the study of mirror symmetry). 
Indeed, this is a key point for the proof of homological mirror 
symmetry for K3 surfaces, for which we refer to \cite{BBS}.

The Maslov class so far constructed does not depend on the choice of 
a canonical projection, from which one could determine the singular 
locus (as usually happens when one considers Lagrangian embedding in 
cotangent bundles over an arbitrary Riemannian manifold). However, it 
is still possible to determine, rather then the singular locus,  
the {\em homology class} $[Z]\in H_{n-1}(\L, Z)$ of a Òsingular locusÓ, 
just considering the Poincare' dual to $[\mu_{\L}]$, and setting 
$[Z]:=Pd([\mu_{\L}])$ ($Pd$ stands for Poncare' duality). We have said 
Òa singular locusÓ, because $Z$ is not determined at all uniquely, 
but only up to its homology class; in spite of this one could take as 
singular locus any representative of $[Z]$. So it makes sense to 
speak of 
a singular locus, even if there is no projection to which to refer it.

It is clear that it is not possible to extend our definition of Maslov 
class for Lagrangian embedding in arbitrary symplectic manifolds; even 
the construction of Fukaya (which is specifically designed for Maslov 
index of closed loops only on BS orbits) needs several assumption 
such that the ambient manifold admits a Òprequantum bundleÓ and so on. 
We are thus tempted to suggest the following alternative description:
we would like to define the Maslov class for a Lagrangian embedding in 
{\em any} symplectic manifold $(X,\o)$, via the mean curvature 
representation $i_{H}\o$. Two problems arise following this approach. 
First of all, to define the mean curvature vector field $H$ it is 
necessary to fix a Riemannian metric on $X$; as it is well known, on 
any symplectic manifold one has lots of Riemannian metrics 
$g_{J}(X,Y):=\o(X,JY)$, constructed using the given symplectic form 
$\o$ and choosing an $\o$-compatible almost complex structure $J$; 
(recall that the set of $\o$-compatible almost complex structures on a 
given symplectic manifold is always non empty and contractible, see 
\cite{GR}). What is the ÒrightÓ choice for $g_{J}$? 

Once we have fixed the right metric,
the second problem is related to the closure of the 1-form $i_{H}\o$, 
considered as a form on $\L$; 
indeed there is no reason, a priori, for which $i_{H}\o$ has to be 
closed. We are thus led to the following:

{\em Conjecture}: Having fixed the Lagrangian embedding $j:\L 
\hookrightarrow X$, on any symplectic manifold $(X,\o)$ there exists at least 
one Riemannian metric $g_{J}$ built up from an  
$\o$-compatible almost complex structure $J$, such that the 1-form 
$i_{H}\o$, considered as a form on $\L$ is closed. Multiplying the 
corresponding cohomological class $[i_{H}\o]$ for a suitable constant 
in such a way that it is integer valued, we call this class the 
{\em Maslov-Morvan class} of the Lagrangian submanifold $\L$. 
  
  It does not seem possible to give an interpretation of this 
  conjectured Maslov-Morvan class via the universal Maslov class, as we have 
done for Calabi-Yau manifolds, since, in general, we have no control on 
the holonomy of $g_{J}$.

Clearly, the study of the relations between the conjectured class 
$[i_{H}\o]$ and the ordinary Maslov class for a Lagrangian embedding in cotangent 
bundles (via Morse families) deserves further effort and is left for 
future investigations.

\bigskip

{\bf Acknowledgements}: It is a pleasure to thank E. Aldrovandi, 
U.Bruzzo, F. Cardin and M. Spera for encouragement and useful discussions.


\begin{thebibliography}{99}

\bibitem{A1} V.I. Arnol'd, {\it Characteristic classes entering in quantization 
conditions}, Funct. Anal. Appl., {\bf 1}, pag.1-13,1967.

\bibitem{A2} V.I. Arnol'd and A.B. Givental', {\it Symplectic Geometry}, in 
Dynamical Systems 4, Springer Verlag, 1990.

\bibitem{BBS} C. Bartocci, U. Bruzzo and S. Sanguinetti, {\it 
Categorial Mirror Symmetry for K3 Surfaces}, Comm. Math. Phys., in 
press.

\bibitem{Besse} A.L. Besse, {\it Einstein Manifolds}, Modern Surveys 
in Mathematics, Springer Verlag, 1987.

\bibitem{BC} R.L. Bishop, R.J. Crittenden, {\it Geometry of 
Manifolds}, Academic Press, 1964.

\bibitem{BT} R. Bott, L.W. Tu, {\it Differential forms in Algebraic 
Topology}, Springer Verlag, 1982.

\bibitem{F} K. Fukaya, {\it Morse homotopy, A$^{\infty}$-category and 
Floer homologies}, preprint MSRI, Berkeley, 1993.

\bibitem{GR} M. Gromov, {\it Pseudo holomorphic curves in symplectic
 manifolds}, Invent. Math. 82 (1985), 307-347. 

\bibitem{GS} V. Guillemin and S. Sternberg, {\it Geometric 
Asymptotics}, Mathematical surveys, n. 14, AMS.

\bibitem{HL} R. Harvey, H.B. Lawson, {\it Calibrated Geometries}, 
Acta Math., 148 (1982), 47-157.

\bibitem{KM} M.V. Karashev, V.P. Maslov, {\it Nonlinear Poisson 
brackets, Geometry and Quantization}, AMS, 1993.

\bibitem{K} M. Kontsevich, {\it Homological Algebra of Mirror 
Symmetry}, alg-geom/9411018.

\bibitem{M} V.P. Maslov, {\it Perturbation Theory and Asymptotic Methods}, 
Dunod, Paris,1972.

\bibitem{Mor} J-M. Morvan, {\it Classe de Maslov d' une immersion 
lagrangienne et minimalitŽ}, C.R. Acad. Sc. Paris, t.292, SŽrie I,1981.

\bibitem{SYZ} A. Strominger, S.-T. Yau, E. Zaslow, {\it Mirror 
Symmetry is T-duality}, Nucl. Phys. B479, (1996), 243-259.

\bibitem{Trof} V.V. Trofimov, {\it Generalized Maslov classes on the 
path space of a symplectic manifold}, Proceedings of the Steklov 
Institute of Math., 1995, Issue 4.

\bibitem{V} C. Voisin, {\it SymŽtrie Miroir}, SMF, Panoramas et 
Synthses, 1996.

\bibitem{W} A. Weinstein, {\it Lectures on Symplectic Manifolds}, 
Regional conference series in mathematics, A.M.S., Providence Rhode 
Island, 1977.

\bibitem{YK} K. Yano, M. Kon, {\it Anti-invariant submanifolds}, 
Lecture Notes in Pure and Applied Math., Vol. 21, Marcel Dekker Inc., 
1976.
\end{thebibliography}
\end{document}